\documentclass[a4paper,12pt,reqno]{amsart}

\usepackage{amsmath}
\usepackage{amssymb}
\usepackage{amsfonts}
\usepackage{graphicx}
\usepackage[colorlinks]{hyperref}
\renewcommand\eqref[1]{(\ref{#1})} 

\graphicspath{ {images/} }
\title[Global existence and blow-up of solutions]{Global existence and blow-up of solutions to pseudo-parabolic equation for Baouendi-Grushin operator}

\author[A. Dukenbayeva]{Aishabibi Dukenbayeva}
\address{
  Aishabibi Dukenbayeva:
    \endgraf
  Suleyman Demirel University
  \endgraf
  Kaskelen, Kazakhstan
  \endgraf
  and
  \endgraf
  Institute of Mathematics and Mathematical Modeling
  \endgraf
 Almaty, Kazakhstan
  \endgraf
 {\it E-mail address} {\rm aishabibi.dukenbayeva@gmail.com}
  }


\thanks{This research is funded by the Committee of Science of the Ministry of Science and Higher Education of the Republic of Kazakhstan (Grant No. AP14972714).}

\newtheoremstyle{theorem}
{10pt}          
{10pt}  
{\sl}  
{\parindent}     
{\bf}  
{. }    
{ }    
{}     
\theoremstyle{theorem}

\numberwithin{equation}{section}
\theoremstyle{plain}
\newtheorem{thm}{Theorem}[section]

\newtheorem{lem}[thm]{Lemma}

\theoremstyle{definition}

\newtheorem{rem}[thm]{Remark}

\newtheoremstyle{defi}
{10pt}          
{10pt}  
{\rm}  
{\parindent}     
{\bf}  
{. }    
{ }    
{}     
\theoremstyle{defi}

\begin{document}
	\begin{abstract}
		In this note, we study a global existence and blow-up of the positive solutions to the initial-boundary value problem of the nonlinear pseudo-parabolic equation for the Baouendi-Grushin operator. The approach is based on the concavity argument and the Poincar\'e inequality related to the Baouendi-Grushin operator from \cite{SY20}, inspired by the recent work \cite{RST23_strat}.
	\end{abstract}
	\maketitle
	
	\section{Introduction}
Consider the initial-boundary problem of the nonlinear pseudo-parabolic equation
	\begin{align}\label{main_eqn_p_Euc}
	\begin{cases}
	u_t - \Delta u_t -\Delta u = f(u), \,\,\, & (x,y) \in \Omega,\,\, t>0, \\ 
	u(x,y,t)  =0,  \,\,\,& (x,y)\in \partial \Omega,\,\, t>0, \\
	u(x,y,0)  = u_0(x,y)\geq 0,\,\,\, & (x,y) \in \overline{\Omega},
	\end{cases}
	\end{align}
	where $f$ is locally Lipschitz continuous on $\mathbb{R}$, $f(0)=0$, and such that $f(u)>0$ for $u>0$, where $\Omega$ is a smoothly bounded domain in $\mathbb{R}^n$. Here, $u_0$ is a non-negative and non-trivial function in $C^1(\overline{\Omega})$ with $u_0(x,y)=0$ on the boundary $\partial \Omega$. 

 	The energy of isotropic materials can be described using a pseudo-parabolic equation \cite{CG68}. Pseudo-parabolic equations are also used to model certain wave processes \cite{BBM72}, as well as the filtration of two-phase flow in porous media considering dynamic capillary pressure \cite{Baren}. Numerous researchers have investigated the global existence and finite-time blow-up of solutions to pseudo-parabolic equations in both bounded and unbounded domains, as referenced in works such as \cite{Korpusov1, Korpusov2, Long, Luo, Peng, Xu1, Xu2, Xu3} and others.

  In this note, we consider an extension of the problem \eqref{main_eqn_p_Euc} to the Baouendi-Grushin setting. 
 	
		Let $z:=(x,y):=(x_{1},...,x_{m}, y_{1},...,y_{k})\in \mathbb{R}^{m}\times\mathbb{R}^{k}$ with $m,k\geq1$ and $m+k=n$. In this setting, we define the corresponding sub-elliptic gradient on $\mathbb{R}^{m+k}$ by 
  \begin{equation}\label{subgrad}
\nabla_{\gamma}:=(X_{1},...,X_{m}, Y_{1},...,Y_{k})=(\nabla_{x}, |x|^{\gamma}\nabla_{y})
\end{equation}
and Baouendi-Grushin operator by
  \begin{equation}\label{Grush_op}
\Delta_{\gamma}=\sum_{i=1}^{m}X_{i}^{2}+\sum_{j=1}^{k}Y_{j}^{2}=\Delta_{x}+|x|^{2\gamma}\Delta_{y}=\nabla_{\gamma}\cdot \nabla_{\gamma},
\end{equation}
  where
$$X_{i}=\frac{\partial}{\partial x_{i}}, \;i=1,...,m, \;\;\; Y_{j}=|x|^{\gamma}\frac{\partial}{\partial y_{j}},\;\gamma\geq0,\;j=1,...,k.$$

Now we are ready to state our problem in the Baouendi-Grushin setting: Let $D \subset \mathbb{R}^{m+k}$ be a bounded domain (open and connected) supporting the divergence formula and $D \backslash \{(x, y) \in \overline{D}: x=0\}$ consists of only one connected component. We consider the problem
	\begin{align}\label{main_eqn_p}
	\begin{cases}
	u_t - \Delta_{\gamma} u_t -\Delta_{\gamma} u = f(u), \,\,\, & (x,y) \in D,\,\, t>0, \\ 
	u(x,y,t)  =0,  \,\,\,& (x,y)\in \partial D,\,\, t>0, \\
	u(x,y,0)  = u_0(x,y)\geq 0,\,\,\, & (x,y) \in \overline{D},
	\end{cases}
	\end{align}
	where $f$ is locally Lipschitz continuous on $\mathbb{R}$, $f(0)=0$, and such that $f(u)>0$ for $u>0$. As above, we assume that $u_0$ is a non-negative and non-trivial function in $C^1(\overline{D})$ with $u_0(x,y)=0$ on the boundary $\partial D$. In this note, we are interested in global existence and blow-up of the positive solutions to \eqref{main_eqn_p}.

Veron and Pohozaev in \cite{PohVer} studied blow-up phenomena for the following equation on the Heisenberg group $\mathbb{H}^{n}$:
	$$
	\frac{\partial u(x,t)}{\partial t}-\mathcal{L} u(x,t)=|u(x,t)|^{p},\,\,\,\,\,(x,t)\in \mathbb{H}^{n} \times(0,+\infty),
	$$
 where $\mathcal{L}$ is the sub-Laplacian on $\mathbb{H}^{n}$. We also refer to \cite{AAK1, AAK2, DL1, JKS1, JKS2} for blow-up type results for semi-linear diffusion and pseudo-parabolic equations on the Heisenberg group as well as to \cite{RY} for the Fujita exponent on general unimodular Lie groups. Our research is inspired by the work \cite{RST23_strat} where the authors obtained global existence and blow-up type results for the problem \eqref{main_eqn_p} with a sub-Laplacian on stratified Lie groups.

It is worth mentioning that the Baouendi-Grushin operator for an even positive integer $\gamma$ is a sum of squares of smooth vector fields satisfying the H\"{o}rmander rank condition
$${\rm rank} \;{\rm Lie} [X_{1},...,X_{m}, Y_{1},...,Y_{k}]=n.$$
So, from \eqref{Grush_op} one can note that unlike the sub-Laplacian on the Heisenberg group or unimodular Lie groups, the Baouendi-Grushin operator not always can be represented by sum of squares of H\"{o}rmander's vector fields. 
Recall also that the anisotropic dilation attached to the Baouendi-Grushin operator is defined by 
$$\delta_{\lambda}(z)=(\lambda x, \lambda^{1+\gamma} y), \quad \lambda>0,$$
and the homogeneous dimension with respect to $\delta_{\lambda}$ is defined by
\begin{equation}\label{hom_dim}
Q=m+(1+\gamma)k.
\end{equation}
A change of variables formula for the Lebesgue measure implies that
$$d\circ\delta_{\lambda}(x,y)=\lambda^{Q}dxdy.$$

Let $H_{0}^{1, \gamma}(D)$ the Sobolev space obtained as completion of $C_{0}^{\infty}(D)$ with respect to the norm
$$
\|f\|_{H_{0}^{1, \gamma}(D)}:=\left(\int_{D}\left|\nabla_{\gamma} f\right|^{2} d x d y\right)^{\frac{1}{2}}.
$$

Thus, the first result of this note on the blow-up property takes the form:
	\begin{thm}\label{thm_p>21}
				Assume that
		\begin{equation}\label{hyp-blow-p}
		\alpha F(u) \leq u f(u) + \beta u^{2} +\alpha\theta,\,\,\, u>0,
		\end{equation}
		where 
		\begin{equation*}
		F(u) = \int_{0}^u f(s)ds,
		\end{equation*}
		for some
		\begin{align}\label{(1.7)}
		\alpha>2 \,\, &\text{ and } \,\,	0<\beta\leq \lambda_{1}\frac{(\alpha - 2)}{2} \\
		\theta>0, \nonumber
		\end{align}
and $\lambda_{1}$ is the first eigenvalue of the Dirichlet Baouendi-Grushin operator on $D$.
		Assume also that the initial data $u_0 \in L^{\infty}(D)\cap H_{0}^{1, \gamma}(D)$ satisfies
		\begin{equation}\label{F(0)}
		\mathcal{F}_0:= -\frac{1}{2} \int_{D} |\nabla_{\gamma} u_0(x,y)|^2 dxdy + \int_{D} (F(u_0(x,y))-\theta) dxdy >0.
		\end{equation}
		Then any positive solution $u$ of \eqref{main_eqn_p} blows up in finite time $T^*,$ that is, there exists 
		\begin{equation}
		0<T^*\leq \frac{M}{\sigma \int_{D} u_0^2 +  |\nabla_{\gamma} u_0|^2 dxdy},
		\end{equation}
		such that
		\begin{equation}
		\lim_{t\rightarrow T^*}\int_{0}^t \int_{D} [u^{2} + |\nabla_{\gamma} u|^2 ]dxdy d\tau = +\infty,
		\end{equation}
		where $\sigma=\sqrt{\frac{\alpha}{2}}-1>0$ and 
		\begin{equation*}
		M = \frac{(1+\sigma)\left( 1+ \frac{1}{\sigma}\right)\left( \int_{D} u^{2}_0 + |\nabla_{\gamma} u_0|^2 dxdy \right)^2}{2\alpha \mathcal{F}_0}.
		\end{equation*}
	\end{thm}
 \begin{rem} The authors in \cite{Chung-Choi} used the same condition on $f(u)$ for a parabolic equation with the classical Laplacian (see also \cite{PhP-06} and \cite{Band-Brun} for particular cases). We can also refer to \cite{RST23_strat}, \cite{ST-21}  and \cite{Bolys_wave} for recent papers on such conditions.
	\end{rem}
 \begin{rem}
				    For more details on spectral properties of the Dirichlet Baouendi-Grushin operator we refer to e.g. \cite{AHKP08},  \cite{MSV15} and \cite{MP09}.
				\end{rem}
The following result indicates that, for certain nonlinearities, positive solutions can be controlled when they exist.  
	\begin{thm}\label{thm_GEp1}
				Assume that 
		\begin{equation}\label{G1}
		\alpha F(u) \geq u f(u) + \beta u^{2} +\alpha\theta, \,\,\, u>0,
		\end{equation}
		where 
		\begin{equation*}
		F(u) = \int_{0}^u f(s)ds,
		\end{equation*}
		for some
		\begin{equation}
		\beta \geq\frac{2-\alpha}{2} \,\,\, \text{ and } \,\,\,\alpha \leq 0, \,\, \theta \geq 0.
		\end{equation} 
		Let the initial data $u_0 \in L^{\infty}(D)\cap H_{0}^{1, \gamma}(D)$ satisfy 
		\begin{equation}
		\mathcal{F}_0:= -\frac{1}{2} \int_{D} |\nabla_{\gamma} u_0(x,y)|^2 dxdy +\int_{D} (F(u_0(x,y))-\theta) dxdy >0. 
		\end{equation}
		If $u$ is a positive local solution of the problem \eqref{main_eqn_p}, then it is global with the property
		\begin{equation*}
		\int_{D}( u^{2} + |\nabla_{\gamma} u|^2) dxdy \leq \exp({-(2-\alpha)t})\int_{D}( u^{2}_0 + |\nabla_{\gamma} u_0|^2) dxdy.
		\end{equation*}
	\end{thm}
The following Poincar\'e inequality established in \cite{SY20} plays an important role for our analysis:
\begin{lem}\label{lem1}
			 			 Let $D \subset \mathbb{R}^{m+k}$ be a bounded domain (open and connected) supporting the divergence formula and $D \backslash \{(x, y) \in \overline{D}: x=0\}$ consists of only one connected component. For every function $u \in H_{0}^{1, \gamma}(D)$ we have 
			\begin{equation}
			\int_{D} |\nabla_{\gamma} u|^2 dxdy \geq \lambda_{1} \int_{D} |u|^2 dxdy,
			\end{equation}
			where $\lambda_{1}$ is the first eigenvalue of the Dirichlet Baouendi-Grushin operator on $D$.
				\end{lem}

	.

	\section{Proofs}\label{sec2}
	Let us begin with the proof of Theorem \ref{thm_p>21} on the blow-up property of the problem \eqref{main_eqn_p}.

	\begin{proof}[Proof of Theorem \ref{thm_p>21}]
Denote
		\begin{align}\label{eq-E_itself}
		E(t) &:= \int_{0}^t \int_{D}( u^{2} + |\nabla_{\gamma} u|^2) dxdy d\tau + M, \,\, t\geq 0,
		\end{align}
		where $M$ is a positive constant to be chosen later. A direct calculation implies that
		\begin{align*}
		E'(t) = \int_{D}( u^2 +  |\nabla_{\gamma} u|^2 )dxdy = \int_{0}^t \frac{d}{d\tau}\int_{D} (u^2 + |\nabla_{\gamma} u|^2)dxdy d\tau + \int_{D} (u^{2}_0 + |\nabla_{\gamma} u_0|^2) dxdy. 
		\end{align*}
  It implies for arbitrary $\delta>0$ that 
		\begin{multline}\label{eq-E_first_der}
		(E'(t))^2\leq  (1+\delta)\left( \int_{0}^t \frac{d}{d\tau}\int_{D}(u^2 + |\nabla_{\gamma} u|^2) dxdy d\tau \right)^2\\
		 + \left( 1+ \frac{1}{\delta}\right)\left( \int_{D} (u^{2}_0 + |\nabla_{\gamma} u_0|^2) dxdy \right)^2.
		\end{multline}
By making use of \eqref{hyp-blow-p} and integration by parts, we obtain for $E''(t)$ that
		\begin{align*}
		E''(t) &= 2\int_{D} u u_t dxdy + \int_{D} (|\nabla_{\gamma} u|^2)_t dxdy\\
		& =2 \int_{D}  ( u\Delta_{\gamma}u+  u \nabla_{\gamma} \cdot (\nabla_{\gamma}  u_t) + uf(u))dxdy + \int_{D} (|\nabla_{\gamma} u|^2)_t dxdy\\
		& = -2 \int_{D} (|\nabla_{\gamma} u|^2 + \nabla_{\gamma} u \cdot \nabla_{\gamma} u_t ) dxdy + 2\int_{D} uf(u)dxdy + \int_{D} (|\nabla_{\gamma} u|^2)_t dxdy
		\\
		& \geq - 2 \int_{D} |\nabla_{\gamma} u|^2 dxdy + 2 \int_{D} \left( \alpha F(u) -\beta u^{2}  -\alpha \theta \right) dxdy\\
		& =2 \alpha \left( -\frac{1}{2} \int_{D} |\nabla_{\gamma} u|^p dxdy + \int_{D} (F(u)-\theta) dxdy \right)  \\
		&+ \frac{2(\alpha-2)}{2} \int_{D} |\nabla_{\gamma} u|^2 dxdy -2\beta \int_{D} u^{2} dxdy
  \\ & \geq 2\alpha \left( -\frac{1}{2} \int_{D} |\nabla_{\gamma} u|^2 dxdy + \int_{D} (F(u)-\theta) dxdy \right)  \\
		& + 2\left(\lambda_{1}\frac{(\alpha -2)}{2} - \beta \right)\int_{D} u^{2} dxdy\\
		&  \geq 2\alpha\left( -\frac{1}{2} \int_{D} |\nabla_{\gamma} u|^2 dxdy + \int_{D} (F(u)-\theta) dxdy \right)=: 2\alpha \mathcal{F}(t),
		\end{align*}
  where we have applied the assumption \eqref{(1.7)} and the Poincar\'e inequality from Lemma \ref{lem1} in the last two estimates above. Note that for $\mathcal{F}_{0}$ from \eqref{F(0)} we have $\mathcal{F}(0)=\mathcal{F}_{0}$. Then, taking into account \begin{equation*}
		\mathcal{F}(t) = \mathcal{F}(0) + \int_{0}^t \frac{d \mathcal{F}(\tau)}{d\tau}d\tau,
		\end{equation*}
  we continue the estimation above as follows
  \begin{align}\label{E_2_deriv}
		E''(t) &=2\alpha \mathcal{F}(t) \nonumber\\&=
  2\alpha\left(\mathcal{F}(0) + \int_{0}^t \frac{d \mathcal{F}(\tau)}{d\tau}d\tau\right) \nonumber\\&= 2\alpha\left(\mathcal{F}(0) - \frac{1}{2} \int_{0}^t\int_{D} \frac{d}{d\tau}|\nabla_{\gamma} u|^2 dxdy d\tau  + \int_{0}^t \int_{D} \frac{d}{d\tau} (F(u)-\theta) dxdy d\tau\right) \nonumber\\
		& =  2\alpha\left(\mathcal{F}(0) -  \int_{0}^t\int_{D} |\nabla u \cdot \nabla_{\gamma} u_{\tau} dxdy d\tau +  \int_{0}^t \int_{D} F_u(u) u_{\tau} dxdy d\tau\right) \nonumber\\
		& = 2\alpha\left(\mathcal{F}(0) +  \int_{0}^t \int_{D} [\Delta_{\gamma} u + f(u)]u_{\tau} dxdy d\tau\right) \nonumber\\
		& = 2\alpha\left(\mathcal{F}(0) +  \int_{0}^t \int_{D} u_{\tau}^2 - u_{\tau} \nabla_{\gamma} \cdot ( \nabla_{\gamma} u_{\tau})dxdy d\tau\right) \nonumber\\
		& = 2\alpha\left(\mathcal{F}(0) + \int_{0}^t \int_{D}  u_{\tau}^2  + |\nabla_{\gamma} u_{\tau}|^2dxdy d\tau\right).
  \end{align}
  Combining \eqref{eq-E_itself}, \eqref{eq-E_first_der}, \eqref{E_2_deriv} and taking $\sigma=\delta=\sqrt{\frac{\alpha}{2}}-1>0$, we obtain
		\begin{align}\label{main_est1}
		&E''(t) E(t) - (1+\sigma) (E'(t))^2\nonumber\\
		&\geq  2\alpha M\mathcal{F}(0)+ 2\alpha \left(\int_{0}^t \int_{D}  (u_{\tau}^2  + |\nabla_{\gamma} u_{\tau}|^2)dxdy d\tau \right) \left(\int_{0}^t \int_{D} (u^{2} + |\nabla_{\gamma} u|^2) dxdy d\tau \right)\nonumber\\
		&- (1+\sigma) (1+\delta)\left( \int_{0}^t \frac{d}{d\tau}\int_{D}(u^2 + |\nabla_{\gamma} u|^2) dxdy d\tau \right)^2 \nonumber\\&- (1+\sigma) \left( 1+ \frac{1}{\delta}\right)\left( \int_{D} (u^{2}_0 + |\nabla_{\gamma} u_0|^2) dxdy \right)^2 \nonumber\\
		&= 2\alpha M\mathcal{F}(0) - (1+\sigma)\left( 1+ \frac{1}{\delta}\right)\left( \int_{D} (u^{2}_0 + |\nabla_{\gamma} u_0|^2) dxdy \right)^2\nonumber\\
		& + 2\alpha \left( \left(\int_{0}^t \int_{D}  (u_{\tau}^2  + |\nabla_{\gamma} u_{\tau}|^2)dxdy d\tau \right) \left(\int_{0}^t \int_{D} (u^{2} + |\nabla_{\gamma} u|^2) dxdy d\tau \right) \right. \nonumber\\
		&- \left.\left(\int_{0}^t \int_{D} (uu_{\tau} + \nabla_{\gamma} u \cdot \nabla_{\gamma} u_{\tau} )dxdy d\tau \right)^2 \right).
		\end{align}
  Here, let us first show that 
  \begin{align}\label{aux_est2}
		&\left( \int_{0}^t \int_{D}  (u^2  + |\nabla_{\gamma} u|^2)dxdy d\tau \right)\left(\int_{0}^t \int_{D}  (u_\tau^2  + |\nabla_{\gamma} u_\tau|^2)dxdy d\tau \right) \nonumber\\
		&-  \left( \int_{0}^t \int_{D}  (u u_{\tau} + \nabla_{\gamma} u \cdot \nabla_{\gamma} u_{\tau})dxdy d\tau  \right)^2  \geq 0.
		\end{align}
  For this, using H\"older's and Cauchy–Schwarz inequalities we observe that
		\begin{multline*}
		\left( \int_{0}^t\int_{D}  (u u_{\tau} +  \nabla_{\gamma} u\cdot \nabla_{\gamma} u_{\tau})dxdy d\tau  \right)^2 \\
		\leq \left( \int_{D} \left(\int_{0}^t u^2 d\tau\right)^{\frac{1}{2}} \left(\int_{0}^t u_{\tau}^2 d\tau\right)^{\frac{1}{2}}  dxdy\right.   \\\left.+\int_{D} \left(\int_{0}^t |\nabla_{\gamma} u|^2 d\tau\right)^{\frac{1}{2}}\left(\int_{0}^t |\nabla_{\gamma} u_{\tau}|^2 d\tau\right)^{\frac{1}{2}} dxdy  \right)^2\\
		=\left( \int_{D} \left(\int_{0}^t u^2 d\tau\right)^{\frac{1}{2}} \left(\int_{0}^t u_{\tau}^2 d\tau\right)^{\frac{1}{2}}  dxdy\right)^2\\+ \left(\int_{D} \left(\int_{0}^t |\nabla_{\gamma} u|^2 d\tau\right)^{\frac{1}{2}}\left(\int_{0}^t |\nabla_{\gamma} u_{\tau}|^2 d\tau\right)^{\frac{1}{2}} dxdy  \right)^2 
  \end{multline*}
  \begin{multline}\label{aux_est1}
		+ 2\left( \int_{D} \left(\int_{0}^t u^2 d\tau\right)^{\frac{1}{2}} \left(\int_{0}^t u_{\tau}^2 d\tau\right)^{\frac{1}{2}}  dxdy\right)\\\times\left(\int_{D} \left(\int_{0}^t |\nabla_{\gamma} u|^2 d\tau\right)^{\frac{1}{2}}\left(\int_{0}^t |\nabla_{\gamma} u_{\tau}|^2 d\tau\right)^{\frac{1}{2}} dxdy  \right)\\
		\leq \left( \int_{D} \int_{0}^t u^2 d\tau dxdy\right)\left( \int_{D} \int_{0}^t u^2_{\tau} d\tau dxdy\right)\\+ \left( \int_{D} \int_{0}^t |\nabla_{\gamma} u|^2 d\tau dxdy\right)\left( \int_{D} \int_{0}^t |\nabla_{\gamma} u_{\tau}|^2 d\tau dxdy\right)\\
		+ 2\left( \left( \int_{D} \int_{0}^t u^2 d\tau dxdy\right)\left( \int_{D} \int_{0}^t u^2_{\tau} d\tau dxdy\right) \right.\\\left.\times \left( \int_{D} \int_{0}^t |\nabla_{\gamma} u|^2 d\tau dxdy\right)\left( \int_{D} \int_{0}^t |\nabla_{\gamma} u_{\tau}|^2 d\tau dxdy\right)\right)^{\frac{1}{2}}.
		\end{multline}
Then, taking into account \eqref{aux_est1} we can verify \eqref{aux_est2} as follows
  \begin{align*}
		&\left( \int_{0}^t \int_{D}  (u^2  + |\nabla_{\gamma} u|^2)dxdy d\tau \right)\left(\int_{0}^t \int_{D}  (u_\tau^2  + |\nabla_{\gamma} u_\tau|^2)dxdy d\tau \right) \\
		&-  \left( \int_{0}^t \int_{D}  (u u_{\tau} + \nabla_{\gamma} u \cdot \nabla_{\gamma} u_{\tau})dxdy d\tau  \right)^2 \\
		& \geq \left( \left( \int_{D} \int_{0}^t u^2 d\tau dxdy\right)^{\frac{1}{2}}\left( \int_{D} \int_{0}^t |\nabla_{\gamma} u_{\tau}|^2 d\tau dxdy\right)^{\frac{1}{2}} \right. \\
		& \left. - \left( \int_{D} \int_{0}^t |\nabla_{\gamma} u|^2 d\tau dxdy\right)^{\frac{1}{2}}\left( \int_{D} \int_{0}^t  u_{\tau}^2 d\tau dxdy\right)^{\frac{1}{2}} \right)^2 \geq 0,
		\end{align*}
  Thus, using \eqref{aux_est2} in \eqref{main_est1} we arrive at

  \begin{equation*}
		E''(t) E(t) - (1+\sigma) (E'(t))^2\geq 2 \alpha M \mathcal{F}(0) - (1+\sigma)\left( 1+ \frac{1}{\delta}\right)\left( \int_{D} (u^{2}_0 + |\nabla_{\gamma} u_0|^2) dxdy \right)^2,
  \end{equation*}
  which taking into account $\mathcal{F}(0)>0$ implies 
  \begin{equation}\label{aux_est3}
		E''(t) E(t) - (1+\sigma) (E'(t))^2 \geq 0
		\end{equation}
  if we choose $M$ as follows
  \begin{equation*}
		M = \frac{(1+\sigma)\left( 1+ \frac{1}{\delta}\right)\left( \int_{D} (u^{2}_0 + |\nabla_{\gamma} u_0|^2) dxdy \right)^2}{2\alpha \mathcal{F}(0)}.
		\end{equation*}
  From \eqref{aux_est3} one can derive for $t\geq 0$ that
		\begin{equation*}
		\frac{d}{dt} \left( \frac{E'(t)}{E^{\sigma+1}(t)} \right) \geq 0  \Rightarrow 	\begin{cases}
		E'(t) \geq \left( \frac{E'(0)}{E^{\sigma+1}(0)} \right) E^{1+\sigma}(t),\\
		E(0)=M.
		\end{cases}
		\end{equation*}
		Here, for $\sigma =\sqrt{\frac{\alpha}{2}}-1>0$ we conclude that 
		\begin{equation*}
		E(t) \geq \left( \frac{1}{M^{\sigma}}-\frac{ \sigma \int_{D} (u_0^2 +  |\nabla_{\gamma} u_0|^2 )dxdy }{M^{\sigma+1}} t\right)^{-\frac{1}{\sigma}},
		\end{equation*}
		hence the blow-up time $T^*$ satisfies 
		\begin{equation*}
		0<T^*\leq \frac{M}{\sigma \int_{D} (u_0^2 +  |\nabla_{\gamma} u_0|^2) dxdy},
		\end{equation*}
		as desired.		
	\end{proof}
 Let us now prove the global existence result.
	\begin{proof}[Proof of Theorem \ref{thm_GEp1}]
		Denote
		\begin{align*}
		\mathcal{E}(t) :=  \int_{D}( u^{2} + |\nabla_{\gamma} u|^2) dxdy.
		\end{align*}
  Recalling the functional $\mathcal{F}(t)$ from the proof of Theorem \ref{thm_p>21} and using \eqref{G1}, we have
		\begin{align*}
		\mathcal{E}'(t)  &= 2\int_{D} u u_t dxdy + \int_{D} (|\nabla_{\gamma} u|^2)_t dxdy\\
		& =2 \int_{D}  (u\Delta_{\gamma} u + u\nabla_{\gamma} \cdot (\nabla_{\gamma} u_t) + uf(u))dxdy +  \int_{D} (|\nabla_{\gamma} u|^2)_t dxdy\\
		& = -2 \int_{D} (|\nabla_{\gamma} u|^2 + \nabla_{\gamma} u \cdot \nabla_{\gamma} u_t ) dxdy + 2\int_{D} uf(u)dxdy + \int_{D} (|\nabla_{\gamma} u|^2)_t dxdy
		\\
		& \leq 2\alpha \left( -\frac{1}{2} \int_{D} |\nabla_{\gamma} u|^2 dxdy + \int_{D} (F(u)-\theta) dxdy \right)  - \frac{2(2-\alpha)}{2}  \int_{D} |\nabla_{\gamma} u|^2 dxdy \\&-2\beta  \int_{D} u^{2} dxdy \\
		& \leq  2\alpha \left( -\frac{1}{2} \int_{D} |\nabla_{\gamma} u|^2 dxdy + \int_{D} (F(u)-\theta) dxdy \right) \\
		&- (2-\alpha)(E(t) - \int_D u^2 dxdy) -2\beta  \int_{D} u^{2} dxdy,\\
		& =2\alpha \mathcal{F}(t) - (2-\alpha)\mathcal{E}(t)  + (2-\alpha -2\beta)\int_{D}u^2 dxdy.
		\end{align*}
		Using the calculation from \eqref{E_2_deriv} and $\beta \geq \frac{2-\alpha}{2}$ we conclude that
		\begin{align*}
		\mathcal{E}'(t)  + (2-\alpha)\mathcal{E}(t)  \leq 2 \alpha \left( \mathcal{F}_0 +   \int_{0}^t \int_{D}  (u_{\tau}^2  + |\nabla_{\gamma} u_{\tau}|^2)dxdy d\tau \right)\leq 0,
		\end{align*}
		 which means
		\begin{equation*}
		\mathcal{E}(t)  \leq \exp({-(2-\alpha)t})\mathcal{E}(0),
		\end{equation*}
		as desired.
	\end{proof}


\begin{thebibliography}{NZW01}
	

				
	\bibitem{AAK1} B. Ahmad, A. Alsaedi, M. Kirane, M. Al-Yami.
\newblock {\em Nonexistence results for higher order pseudo-parabolic equations in the Heisenberg group.}
\newblock Math. Methods Appl. Sci., 40:1280--1287, 2017.

\bibitem{AAK2} B. Ahmad, A. Alsaedi, M. Kirane. 
\newblock {\em Blow-up of solutions to parabolic inequalities in the Heisenberg group.}
\newblock Electron. J. Differential Equations, 2015: Art. No. 167, 2015.

\bibitem{AHKP08}
C.T. Anh, P.Q. Hung, T.D. Ke, T.T. Phong. 
\newblock {\em Global attractor for a semilinear
parabolic equation involving Grushin operator.}
\newblock Electron. J. Differential Equations, 2008: Art. No. 32, 2008.



\bibitem{Band-Brun}
		C. Bandle, H. Brunner. 
\newblock {\em Blow-up in diffusion equations, a survey.}
\newblock J. Comput. Appl. Math., 97:3--22, 1998.

		\bibitem{Baren} G.I. Barenblatt, J. Garcia-Azorero, A. De Pablo, J.L. Vazquez. 
\newblock {\em Mathematical model of the non-equilibrium water-oil displacement in porous strata.}
\newblock Appl. Anal., 65:19--45, 1997.
		
  \bibitem{BBM72} T.B. Benjamin, J.L. Bona, J.J. Mahony. 
\newblock {\em Model equations for long waves in nonlinear dispersive systems.}
\newblock Philos. Trans. Royal Soc. London. Ser. A, 272:47--78, 1972.


		
	
		\bibitem{Chung-Choi}
		S.-Y. Chung, M.-J. Choi.
\newblock {\em A new condition for the concavity method of blow-up solutions to $p$-Laplacian parabolic equations.}
\newblock J. Differential Equations, 265:6384--6399, 2018.



\bibitem{CG68} 
		P.J. Chen, M.E. Gurtin.
\newblock {\em On a theory of heat conduction involving two temperatures.}
\newblock Z. Angew. Math. Phys., 19:614--627, 1968.




\bibitem{DL1} L. D'Ambrosio.
\newblock {\em Critical degenerate inequalities on the Heisenberg group.}
\newblock Manuscripta Math., 106:519--536, 2001.












		
		\bibitem{JKS1} M. Jleli, M. Kirane, B. Samet.
\newblock {\em Nonexistence results for a class of evolution equations in the Heisenberg group.}
\newblock Fract. Calc. Appl., 18:717--734, 2015.

\bibitem{JKS2} M. Jleli, M. Kirane, B. Samet.
\newblock {\em Nonexistence results for pseudo-parabolic equations in the Heisenberg group.}
\newblock Monatsh. Math., 180:255--270, 2016.

\bibitem{Korpusov1} M.O. Korpusov, A.G. Sveshnikov.
\newblock {\em On the blow-up of solutions of semilinear pseudoparabolic equations with rapidly growing nonlinearities.}
\newblock Zh. Vychisl. Mat. Mat. Fiz., 45:145--155, 2005. (in russian)
		
\bibitem{Korpusov2} M.O. Korpusov, A.G. Sveshnikov.
\newblock {\em On the blow-up in a finite time of solutions of initial-boundary-value problems for pseudoparabolic equations with the pseudo-Laplacian.}
\newblock Zh. Vychisl. Mat. Mat. Fiz., 45:272--286, 2005. (in russian)
		










		

		
		\bibitem{Long}
		Q.F. Long, J.Q. Chen.
\newblock {\em Blow-up phenomena for a nonlinear pseudo-parabolic equation with nonlocal source.}
\newblock Appl. Math. Lett., 74:181--186, 2017.

\bibitem{Luo}
		P. Luo.
\newblock {\em Blow-up phenomena for a pseudo-parabolic equation.}
\newblock Math. Meth. Appl. Sci., 38:2636--2641, 2015.

\bibitem{MP09}
D.D. Monticelli, K.R. Payne.
\newblock {\em Maximum principles for weak solutions of degenerate elliptic equations with a uniformly elliptic direction.}
\newblock J. Differential Equations, 247:1993--2026, 2009.

\bibitem{MSV15}
M. Mihailescu, D. Stancu-Dumitru, C. Varga.
\newblock {\em On the spectrum of a Baouendi–Grushin type operator: an Orlicz–Sobolev space setting approach.}
\newblock Nonlinear Differ. Equ. Appl., 22:1067--1087, 2015.

\bibitem{PhP-06} 
		G.A. Philippin, V. Proytcheva. 
\newblock {\em Some remarks on the asymptotic behaviour of the solutions of a class of parabolic problems.}
\newblock Math. Methods Appl. Sci., 29:297--307, 2006.


		

		
		\bibitem{Peng}
		X.M. Peng, Y.D. Shang, X.X. Zheng.
\newblock {\em Blow-up phenomena for some nonlinear pseudo-parabolic equations.}
\newblock Appl. Math. Lett., 56:17--22, 2016.

\bibitem{RY}
		M. Ruzhansky, N. Yessirkegenov.
\newblock {\em Existence and non-existence of global solutions for semilinear heat equations and inequalities on sub-Riemannian manifolds, and Fujita exponent on unimodular Lie groups.}
\newblock J. Differential Equations, 308:455--473, 2022.

\bibitem{RST23_strat}
		 M. Ruzhansky, B. Sabitbek, B. Torebek. 
\newblock {\em Global existence and blow-up of solutions to porous medium equation and pseudo-parabolic equation, I. Stratified group.}
\newblock Manuscripta Math., 171:377--395, 2023.

\bibitem{ST-21}
		B. Sabitbek, B. Torebek.
\newblock {\em Global existence and blow-up of solutions to the double nonlinear porous medium equation.}
\newblock Discrete Continuous Dyn. Syst. Ser. B., 44:743--767, 2024.

\bibitem{Bolys_wave}
		B. Sabitbek.
\newblock {\em Global existence and nonexistence of solutions for semilinear wave equation with a new condition.}
\newblock Discrete Continuous Dyn. Syst. Ser. B., 43:2637--2657, 2023.





\bibitem{SY20}
		D. Suragan, N. Yessirkegenov. 
\newblock {\em Sharp remainder of the Poincar\'e inequality for Baouendi--Grushin vector fields.}
\newblock Asian-Eur. J. Math., 16: Art. No. 2350041, 2023.

\bibitem{PohVer} L. V\'{e}ron, S.I. Pohozaev.
\newblock {\em Nonexistence results of solutions of semilinear differential inequalities on the Heisenberg group.}
\newblock Manuscripta Math., 102:85--99, 2000.

\bibitem{Xu1} R.Z. Xu, J. Su.
\newblock {\em Global existence and finite time blow-up for a class of semilinear pseudo-parabolic equations.}
\newblock J. Funct. Anal., 264:2732--2763, 2013.

\bibitem{Xu2} R.Z. Xu, X.C. Wang, Y.B. Yang.
\newblock {\em Blowup and blowup time for a class of semilinear pseudo-parabolic equations with high initial energy.}
\newblock Appl. Math. Lett., 83:176--181, 2018.

  \bibitem{Xu3} R.Z. Xu, X.C. Wang.
\newblock {\em Global existence and finite time blowup for a nonlocal semilinear pseudo-parabolic equation.}
\newblock Adv. Nonlinear Anal., 10:261--288, 2021.
		
	\end{thebibliography}
\end{document}